\documentclass[graybox]{svmult}

\usepackage{mathptmx}        
\usepackage{makeidx}         
\usepackage{cite}
\usepackage{graphicx}        
\usepackage{grffile}
\usepackage{epstopdf}        
\usepackage{amsmath}
\usepackage{multicol}        
\usepackage[bottom]{footmisc}

\makeindex             

\begin{document}

\title*{On the efficient parallel computing of long term reliable trajectories
for the Lorenz system}

\author{I. Hristov$^{1,a}$, R. Hristova$^1$, S. Dimova$^1$, P. Armyanov$^1$, N. Shegunov$^1$, \\ I. Puzynin$^2$,  T. Puzynina$^2$, Z. Sharipov$^{2,b}$, Z. Tukhliev$^2$\\
\vspace{0.5cm}
\emph{$^1$ Sofia University, Faculty of Mathematics and Informatics, Bulgaria}\\
\emph{$^2$ JINR, Laboratory of Information Technologies, Dubna, Russia}\\
\vspace{0.5cm}
\emph{E-mails:  $^a$  ivanh@fmi.uni-sofia.bg \hspace{0.3cm} $^b$ zarif@jinr.ru}}
\authorrunning{I. Hristov, R. Hristova et al.}
\titlerunning{On the efficient computing of long term reliable trajectories
for the Lorenz system}

\maketitle


\abstract {In this work we propose an efficient parallelization of  multiple-precision Taylor series method with variable stepsize and fixed order.
For given level of accuracy the optimal variable stepsize determines higher order of the method than in the case of optimal fixed stepsize.
Although the used order of the method is greater then that in the case of fixed stepsize, and hence the computational work per step is greater,
the reduced number of steps gives less overall work. Also the greater order of the method is beneficial in the sense that it increases the parallel
efficiency.
As a model problem we use the paradigmatic Lorenz system.
With 256 CPU cores in Nestum cluster, Sofia, Bulgaria, we succeed to obtain a correct reference solution in the rather long time interval - [0,11000].
To get this solution we performed two large computations: one computation with 4566 decimal digits of precision and 5240-th order  method,
and second computation  for verification - with 4778 decimal digits of precision and 5490-th order  method.}

\section{Introduction}
\label{sec:2}
Multiple precision Taylor series method is an affordable and very efficient numerical
method for integration of some classes of low dimensional dynamical
systems in the case of  high precision demands \cite{Barrio1, Barrio2}. The method
gives a new powerful tool for theoretical investigation of such systems.

A numerical procedure for computing reliable trajectories
of chaotic systems, called Clean Numerical Simulation (CNS), is proposed by Shijun Liao in \cite{Liao1} and applied for different systems \cite{ Liao2, CNS, Li}.
The procedure is based on multiple precision Taylor series method. The main concept for CNS is the critical predictable time $T_c$, which is a kind of  practical Lyapunov time.
$T_c$ is defined as the time for decoupling of two trajectories computed by two different numerical schemes.
The CNS works as follows. An optimal fixed stepsize is chosen. Then estimates of the required order of the  method  $N$ and the required precision
(the number of exact decimal digits $K$ of the floating point numbers) are obtained.
The optimal order $N$ is estimated by computing the $T_c-N$ dependence by means of the  numerical solutions for fixed large enough $K$.
The estimate of $K$ is obtained by computing the $T_c-K$ dependence by means of the numerical solutions for fixed large enough $N$.
This estimate of $K$ is in fact an estimate for the Lyapunov exponent \cite{Wang}.
For given $T_c$ the solution is then computed with the estimated  $N$ and $K$ and after that one more computation with  higher $N$ and $K$ is performed for verification.
The choice of $N$ and  $K$ ensures that the round-off error and the truncation error are of the same order.

When very high precision and very long integration interval are needed, the computational problem can become pretty large.
In this case the parallelization of the Taylor series method is an important task and needs to be carefully developed.
The first parallelization of CNS is reported in \cite{par1} and later improved in \cite{par2}.
A pretty long reference solution for the paradigmatic Lorenz system, namely in the time interval [0,10000], obtained in about 9 days and 5 hours by using the computational resource of 1200 CPU cores, is given in \cite{Liao3}. However, no details of the parallelization process are given in \cite{par1, par2, Liao3}.
In our recent work \cite{Hybrid} we reported in details a simple and efficient hybrid MPI+OpenMP
parallelization of CNS for the Lorenz system and tested it for the same parameters as those in \cite{Liao3}.
The results show very good efficiency and very good parallel performance scalability of our program.

This work can be regarded as a continuation of our previous work \cite{Hybrid}, where fixed stepsize is used.
Here we make a modification of CNS with a variable stepsize and fixed order following the simple approach given in \cite{Jorba}.
Although the used order of the method is greater then that in the case of fixed stepsize, and hence the computational work per step is greater,
the reduced number of steps gives less overall work. Also the greater order of the method is beneficial in the sense that it increases the parallel
efficiency. With 256 CPU cores in Nestum cluster, Sofia, Bulgaria, we succeed to obtain a correct reference solution in [0,11000]
and in this way we improve the results from \cite{Liao3}.
To obtain this solution we performed two large computations: one computation with 4566 decimal digits of precision and 5240-th order  method,
and second computation  for verification - with 4778 decimal digits of precision and 5490-th order  method for verification.
The computations lasted $\approx$ 9 days and 18 hours and $\approx$ 11 days and 7 hours respectively.
Let us note that the improvement of the numerical algorithm does not change the parallelization strategy
from our previous work \cite{Hybrid}, where the parallelization process is explained in more details.
The difference from the previous parallel program is one additional OpenMP single  section with negligible computational work, which computes the optimal step.

It is important to mention that although our test model is the classical Lorenz system,
the proposed parallelization strategy is rather general - it could be applied as well to a large class of chaotic dynamical systems.

\section{Taylor series method and CNS for the Lorenz system}
\label{sec:3}

We consider as a model problem the classical Lorenz system \cite{Lorenz}:
\begin{equation}
\begin{aligned}
\frac{dx}{dt} &= \sigma(y-x)\\
\frac{dy}{dt} &= Rx - y -xz\\
\frac{dz}{dt} &= xy - bz,
\end{aligned}
\end{equation}
where $R=28$, $\sigma=10$, $b=8/3$  are the standard Salztman's parameters.
For these parameters the system is chaotic.
Let us denote with $x_i, y_i, z_i, i=0, ... , N$ the normalized derivatives
(the derivatives divided by \textbf{i!}) of the approximate solution  at the current time $t$.
Then the N-th order Taylor series method for (1) with stepsize $\tau$ is:
\begin{equation}
\begin{aligned}
x(t+\tau) &\approx x_0 + \sum_{i=1}^{N} x_i \tau^i,\\
y(t+\tau) &\approx y_0 + \sum_{i=1}^{N} y_i \tau^i,\\
z(t+\tau) &\approx z_0 + \sum_{i=1}^{N} z_i \tau^i.
\end{aligned}
\end{equation}

The i-th Taylor coefficients (the normalized derivatives)
are computed as follows. From system (1) we have

$$
\begin{aligned}
x_1 &= \sigma(y_0 - x_0),\\
y_1 &= R x_0 - y_0 - x_0 z_0,\\
z_1 &= x_0 y_0 -b z_0.
\end{aligned}
$$

By applying the Leibniz rule for the derivatives of the product of two functions,
we obtain the following recursive procedure for computing  $x_{i+1}, y_{i+1}, z_{i+1}$  for $i=0,..., N-1$:

\begin{equation}
\begin{aligned}
x_{i+1}  &= \frac{1}{i+1} \sigma (y_i - x_i),\\
y_{i+1}  &= \frac{1}{i+1} (Rx_i -y_i -\sum_{j=0}^{i}x_{i-j}z_j),\\
z_{i+1}  &= \frac{1}{i+1} (\sum_{j=0}^{i}x_{i-j}y_j -bz_i).
\end{aligned}
\end{equation}

To compute the {\it{i+1}}-st coefficient in the Taylor series we need all previous coefficients from 0 to {\it{i}}.
In fact, this algorithm for computing the coefficients of the Taylor series
is called automatic differentiation, or sometimes algorithmic differentiation \cite{Moore}.
It is obvious  that we need $O(N^2)$ floating point operations for computing all coefficients. The subsequent evaluation of Taylor series with Horner's rule needs only $O(N)$ operations.

Let us now explain how we choose the stepsize $\tau$. We use a variable stepsize strategy, which makes the method much more robust then
in the fixed stepsize case. We use a simple strategy taken from \cite{Jorba},
which ensures both the convergence of the Taylor series and the minimization of the computational work per unit time.
If we denote the vector of the normalized derivatives of the solution with $\mathbf{X_i}=(x_i,y_i,z_i)$
and take a safety factor 0.993, then the stepsize $\tau$ is determined by the last two terms of the Taylor expansions \cite{Jorba}:

\begin{equation}
\begin{aligned}
\tau=\frac{0.993}{e^2}\min\left\{   {\left(\frac{1}{{\|\mathbf{X_{N-1}} \|}_{\infty}}\right)}^{\frac{1}{N-1}}, {\left(\frac{1}{{\|\mathbf{X_{N}}\|}}_{\infty}\right)}^{\frac{1}{N}}\right\}
\end{aligned}
\end{equation}

\begin{figure}
\begin{center}

\includegraphics[scale=0.45]{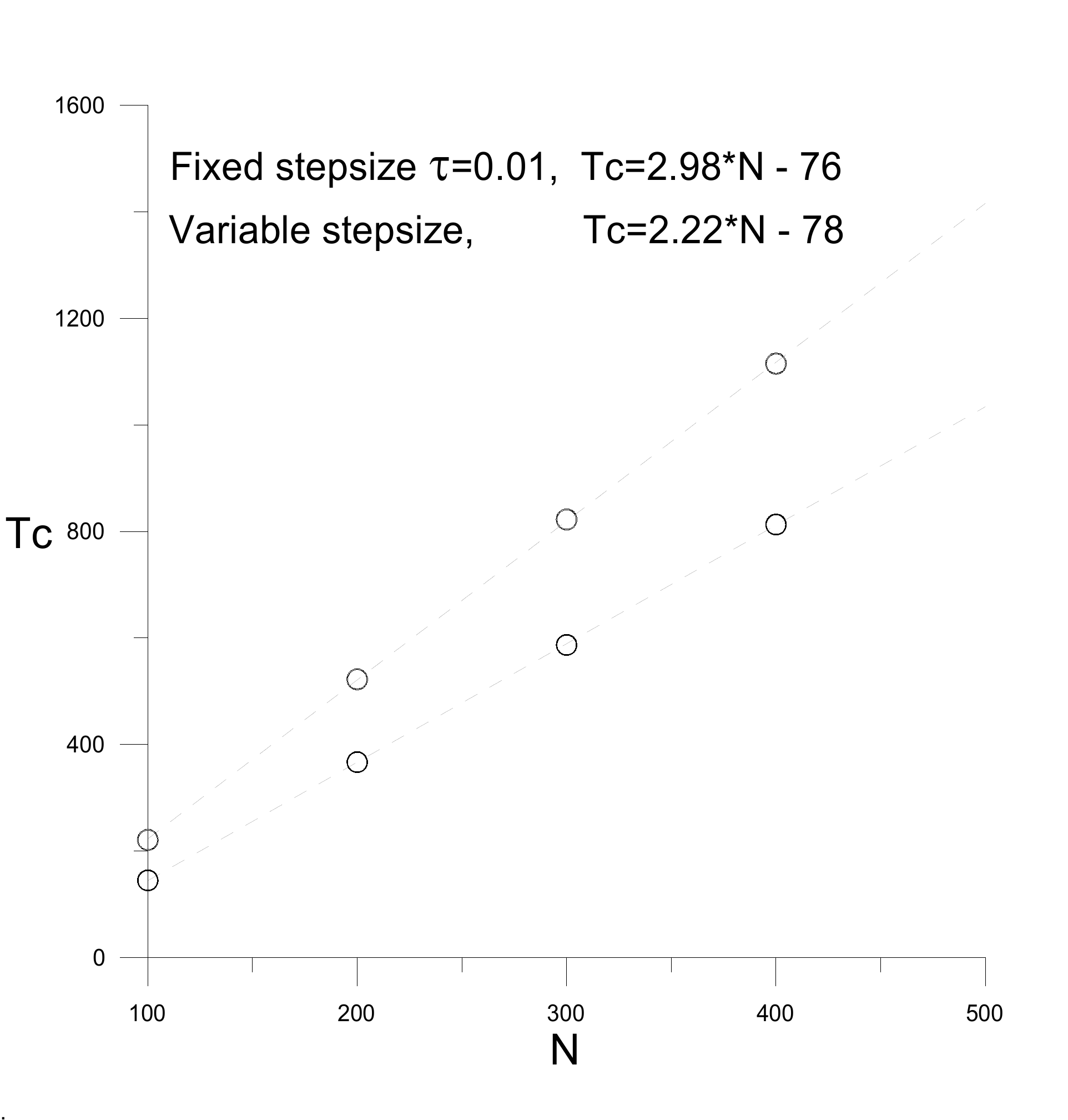}
\caption{$Tc-N$ dependencies for fixed and variable stepsize}
\label{fig:1}
\end{center}
\end{figure}

In \cite{Jorba} the order of the method is determined by the local error tolerance.
However, we do not work explicitly  with some local error tolerance
and also we do not use any explicit dependence between the local and the global error.
Instead of this, as in \cite{Liao1}, we compute an a priori estimate of the needed order of the method for a reliable solution.
As said before, the critical predictable time $T_c$ is defined as the time for decoupling of two trajectories computed by two different numerical schemes
(in this case - by different $N$). The solutions are computed with large enough precision to ensure that the truncation error is the leading one.
As a criterion for decoupling time we choose the time for establishing only 30 correct digits.
The obtained $T_c - N$ dependencies for fixed stepsize $\tau=0.01$ and variable stepsize are shown in Figure 1.
As seen from this figure, the computational work for one step in the case of variable stepsize is $ \approx 80\%$ greater then in the case
of fixed stepsize - ${(2.98/2.22)}^2 \approx 1.80$. However, the reduced number of steps gives less overall work.
Also the greater order of the method is beneficial in the sense that it increases the parallel efficiency.
The reason is that with increasing the order $N$ of the method, the parallelizable part of the work
becomes  relatively even more larger than the serial part and the parallel overhead part.

Similarly, we  compute an a priori estimate of the needed precision by means of computing the $T_c - K$ dependence.
In this case we compare the solutions for different $K$ and large enough $N$.
We obtain the dependence $T_c=2.55K -81$, which is the same, as expected,  for fixed and for variable stepsize.

\section{Parallelization of the algorithm}
The improvement of the numerical algorithm does not change the parallelization strategy
from our previous work \cite{Hybrid}, where the parallelization process is explained in more details.
However, as we will see, the variable stepsize not only decreases the computational work for a given accuracy,
but also gives a higher parallel efficiency.

Let us store the Taylor coefficients in the arrays \textbf{x}, \textbf{y}, \textbf{z} of lengths N+1.
The values of $x_i$ are stored in \textbf{x[i]},
those of $y_i$ in \textbf{y[i]} and those of $z_i$ in \textbf{z[i]}.
As explained in \cite{par1, par2}, the crucial decision for parallelization is
to make a parallel reduction for the two sums in (3). However, in order to reduce the remaining serial part of the code
and hence to improve the parallel speedup from the Amdal's law, we should utilize some limited, but important parallelism.
We compute  \textbf{x[i+1]}, \textbf{y[i+1]}, \textbf{z[i+1]} in parallel.
Moreover, we compute \textbf{x[i+1]} in advance, before computing the sums in (3), when during the reduction
process some of the computational resource is free.
In the same way we compute in advance \textbf{Rx[i]-y[i]} from the formula for \textbf{y[i+1]} and  \textbf{bz[i]} from the formula for \textbf{z[i+1]}.
These computations taken in advance matter, because multiplication is much more expensive then the other used operations (division by an integer number is not so expensive).
The three evaluations by Horner's rule  for the new \textbf{x[0]}, \textbf{y[0]}, \textbf{z[0]} are also done in parallel.

In this work we consider a hybrid MPI+OpenMP strategy \cite{MPI, OpenMP}, i.e. every MPI process creates a team of OpenMP threads.
For multiple precision floating point arithmetic we use GMP library (GNU Multiple Precision library)  \cite{gnu}.
The main reason to consider a hybrid strategy, rather than a pure MPI one, is that OpenMP performs slightly better than MPI on one computational node.
For packaging and unpackaging of the GMP multiple precision types for the MPI messages, we rely on the tiny  MPIGMP library of Tomonori Kouya \cite{Kouya0, Kouya1, Kouya2, Kouya3}.

\begin{figure}
\begin{center}

\includegraphics[scale=1.0]{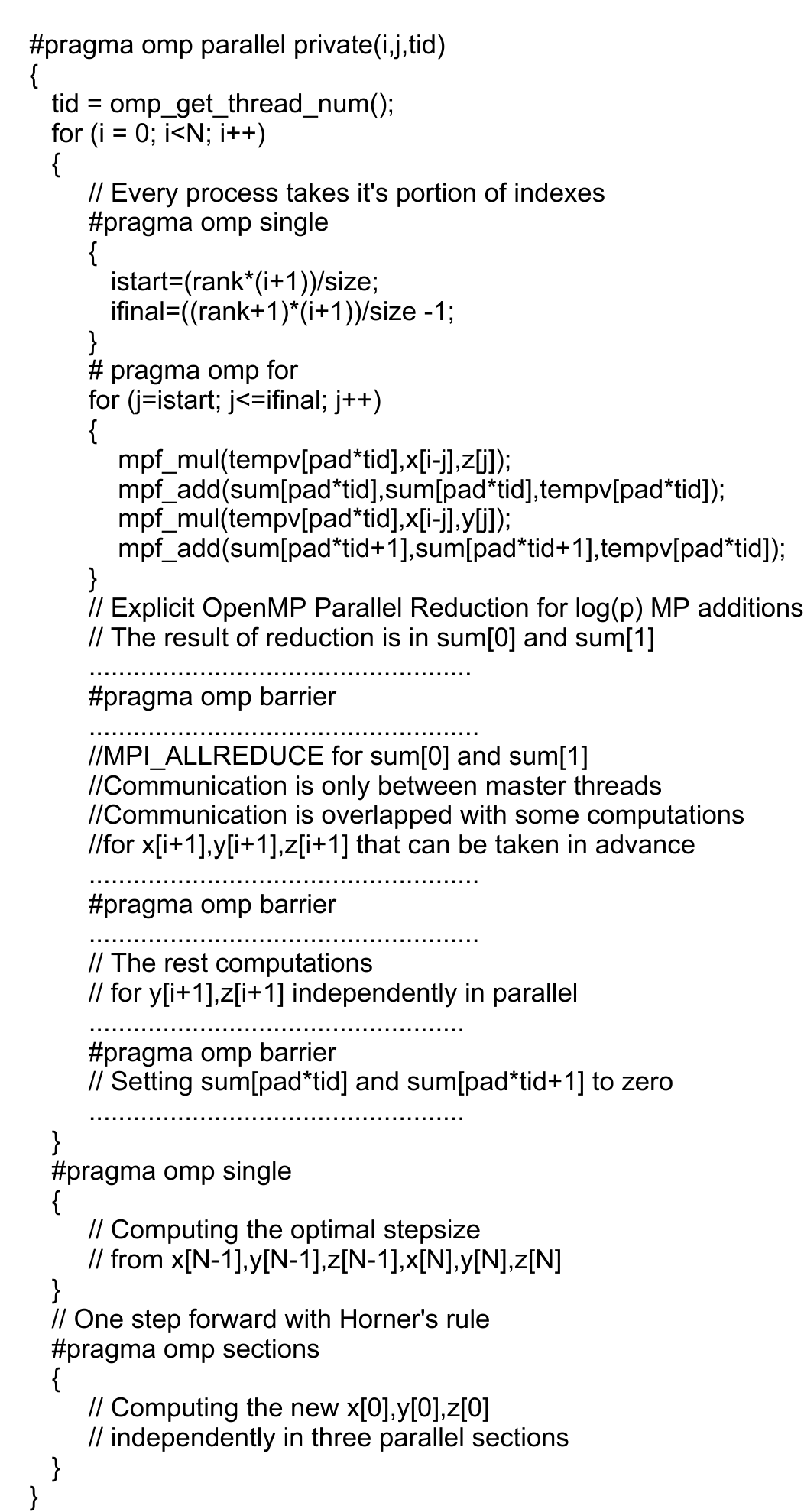}
\caption{The sketch of hybrid MPI+OpenMP code in terms of GMP library.}
\label{fig:1}
\end{center}
\end{figure}

It is important to note that for our problem the pure OpenMP parallelization has its own importance.
First, the programming with OpenMP is easier because it avoids the usage of libraries like MPIGMP.
Second, since the algorithm does not allow domain decomposition,
the  memory needed for one computational node is multiplied by the number of the MPI processes per that node,
while OpenMP needs only one copy of the computational domain and thus some memory is saved.

The sketch of our parallel program is given in Figure 2.
Every thread gets its \textbf{id} and stores it in \textbf{tid} and
then the  loop with index \textbf{i} is performed. Every MPI process takes its portion - the first and the last index controlled by the process.
After that the directive \textbf{\#pragma omp for} shares the work for
the loop between threads.

Although OpenMP has a build-in reduction clause,  we can not use it,
because we use user-defined types for multiple precisions number and user-defined operations.
A manual reduction by applying a standard tree based parallel reduction is done.
We use containers for the partial sums of every thread and these containers are shared.
The containers are stored in the array \textbf{sum}.
We have in addition an array of temporary variables \textbf{tempv} for storing the intermediate results of the multiplications.
To avoid false sharing, a padding strategy is applied \cite{OpenMP}.
At the point where each process has computed its partial sums, we perform  MPI\_ALLREDUCE between the master threads\cite{MPI}.
It is useful to regard MPI\_ALLREDUCE as a continuation of the tree based reduction process, which starts with the OpenMP reduction.
Communications between master threads are overlapped with some computations for \textbf{x[i+1]}, \textbf{y[i+1]}, \textbf{z[i+1]} that can be taken before the computation of the sums in (3) is finished. When the MPI\_ALLREDUCE  is finished, we compute in parallel the remaining operations for \textbf{x[i+1]}, \textbf{y[i+1]}, \textbf{z[i+1]}.

In between the block which computes the Taylor coefficients and the block which computes the new values of \textbf{x[0]}, \textbf{y[0]}, \textbf{z[0]}
in parallel, we compute the new optimal stepsize within an \textbf{omp single} section. While the block for computing
the Taylor coefficients is $O(N^2)$ and the block for evaluations of the polynomials is $O(N)$, this block is only $O(1)$ and hence the work is negligible.
Let us note that the GMP library does not offer a \textbf{power} function for
the computations from formula (4).
The good thing is that we do not need to compute the stepsize with multiple precision and double precision is enough.
So we use the C standard library function \textbf{pow} in double precision.
We do a normalization of the large GMP floating point numbers in order to work in the range of the standard double precision numbers.
The C-code in terms of GMP library of our hybrid MPI+OpenMP program can be downloaded from \cite{radahpc}.

Let us mention that if one half of the OpenMP threads computes one of the sums
in (3) and the other half computes the other sum,
one could also expect some small performance benefit, because for the small indexes \textbf{i} the unused threads
will be less and also the difference from the perfect load balance between threads will be less.
However, the last approach is not general because it strongly depends
on the number of sums for reduction (two in the particular case of the Lorenz system) and the number of available threads.

\section{Computational resources. Performance and numerical results.}
The preparation of the parallel program  and the many tests are performed
in the {\textbf{Nestum} Cluster, Sofia, Bulgaria \cite{nestum} and in the \textbf{HybriLIT} Heterogeneous Platform at the Laboratory of IT of JINR, Dubna, Russia \cite{HybriLIT}.
The large computations for the reference solution in the time interval [0,11000] and the presented results for the performance are from Nestum Cluster.
Nestum is a homogeneous HPC cluster based on two socket nodes. Each node consists of 2 x Intel(R) Xeon(R) Processor E5-2698v3 (Haswell-based processors)  with 32 cores at 2.3 GHz.
We have used  Intel C++ compiler version 17.0, GMP library version 6.2.0, OpenMPI version 3.1.2 and compiler optimization options
-O3 -xhost.

We use the same initial conditions as those in  \cite{Liao3}, namely  $x(0)=-15.8$, $y(0)=-17.48$, $z(0)=35.64$, in order to compare
with the benchmark table in  \cite{Liao3}. We computed a reference solution in the rather long time interval [0,11000]
and repeated the benchmark table up to time 10000. Computing this table by two different stepsize strategies,
is a good demonstration that Clean Numerical Simulation (CNS) is a correct and valuable approach for
computing reliable trajectories of chaotic systems.

We performed two large computations with 256 CPU cores (8 nodes in Nestum).
The first computation is with 4566 decimal digits of precision and 5240-th order  method
(5\% reserve from the a priori estimates). The second computation  is for verification - with 4778 decimal digits of precision and 5490-th order
method (10\% reserve from the a priori estimates).
The first computation lasted $\approx$ 9 days and 18 hours and the second $\approx$ 11 days and 7 hours.
The overall speedup with 256 cores for the first computation is 162.8, for the second - 164.6.

By estimating the time needed for the same accuracy and with fixed stepsize 0.01,
we conclude that by applying variable stepsize strategy we have \textbf{2.1x} speedup.
There are two reasons for this speedup - less overall work and increased parallel efficiency.
Although the work per step in the case of variable stepsize  increases by $ \approx 80\%$, the average stepsize
is $\approx$ 0.034 and thus the overall work is $ \approx 53\%$ from the work in the case of fixed stepsize 0.01.
Also the parallel efficiency increases from 55.5\% up to 63.6\% for the first computation and from 56.2\% up to 64.3\% for the second.
This is because by increasing the order of the method $N$, we increase the amount of the  parallel work,
which mitigates the impact of the serial work and the parallel overhead work.

As we compute the reference solution with some reserve of the estimated $N$ and $K$, we actually obtain
the solution with some more correct digits. The reference solution with 60 correct digits at every 100 time units
can be seen in \cite{radahpc}. The reference solution at $t=11000$ is:\\

x= 6.10629269055689971917782003095370055267185885053970862735508\\

y=-3.33795350928712428173974978144552360814210542698512462640748\\

z= 34.1603471532583648867450334710712261840913307358242610005285

\section{Conclusions}
Parallelized version of multiple precision Taylor series method and particularly the Clean Numerical Simulation
should be used with a variable stepsize strategy as a better alternative of the fixed stepsize one.
An important observation is that variable stepsize not only decreases the computational work for a given accuracy,
but also gives a higher parallel efficiency.

\begin{acknowledgement}
We thank for the opportunity to use the computational resources of the Nestum cluster, Sofia, Bulgaria.
We would like to give our special thanks to Dr. Stoyan Pisov for his great help in using the Nestum cluster
and Prof. Emanouil Atanassov from IICT, BAS for valuable discussions and important remarks on the parallelization process.
We also thank the Laboratory of Information Technologies of JINR, Dubna, Russia for the opportunity to use the computational resources of the HybriLIT Heterogeneous Platform.
The work is supported by a grant of the Plenipotentiary Representative of the Republic of Bulgaria at JINR, Dubna, Russia.
\end{acknowledgement}

\end{document}